\documentclass[letterpaper, 10 pt, conference]{ieeeconf}  

\usepackage{xcolor}
\usepackage{cite}
\usepackage{amsmath,amssymb,amsfonts}
\usepackage{graphicx}
\usepackage{textcomp}
\usepackage[utf8]{inputenc}
\usepackage{float}
\usepackage{subcaption}
\usepackage[ruled, vlined, linesnumbered]{algorithm2e}
\usepackage{booktabs}
\DontPrintSemicolon
\usepackage[left=2cm,right=2cm, top=2cm, bottom=2cm]{geometry}
\usepackage{pgfplots}
\DeclareUnicodeCharacter{2212}{−}
\usepgfplotslibrary{groupplots,dateplot}
\usetikzlibrary{patterns,shapes.arrows}
\pgfplotsset{compat=newest}
\usepackage{url}
\usepackage{courier}
\usepackage{pgf}

\IEEEoverridecommandlockouts
\title{\LARGE \bf
An Almost Feasible Sequential Linear Programming Algorithm
}

\author{David Kiessling$^{1}$, Charlie Vanaret$^3$, Alejandro Astudillo$^{1}$, Wilm Decr\'e$^{1}$, Jan Swevers$^{1}$
\thanks{$^{1}$MECO~Research~Team,~Dept.~of~Mechanical~Engineering,~KU~Leuven and Flanders Make@KU Leuven, 3001 Leuven, Belgium.
        {\tt\small \{david.kiessling, alejandro.astudillovigoya, jan.swevers\}@kuleuven.be}}%
\thanks{$^{3}$Mathematical Algorithmic Intelligence Division, Zuse-Institut Berlin, Berlin, Germany.
{\tt\small vanaret@zib.de}}%
\thanks{This work has been carried out within the framework of the Flanders
Make SBO project DIRAC: Deterministic and Inexpensive Realizations of Advanced
Control.}}

\begin{document}

\maketitle
\thispagestyle{empty}
\pagestyle{empty}

\begin{abstract}
This paper proposes an almost feasible Sequential Linear Programming (afSLP) algorithm. In the first part, the practical limitations of previously proposed Feasible Sequential Linear Programming (FSLP) methods are discussed along with illustrative examples. Then, we present a generalization of FSLP based on a tolerance-tube method that addresses the shortcomings of FSLP. The proposed algorithm afSLP consists of two phases. Phase I starts from random infeasible points and iterates towards a relaxation of the feasible set. Once the tolerance-tube around the feasible set is reached, phase II is started and all future iterates are kept within the tolerance-tube. The novel method includes enhancements to the originally proposed tolerance-tube method that are necessary for global convergence.
afSLP is shown to outperform FSLP and the state-of-the-art solver IPOPT on a SCARA robot optimization problem.
\end{abstract}



\section{Introduction}

In nonlinear model predictive control, the solution of a nonlinear optimization problem must be computed in every sampling period so that the system is given an optimized control signal. Issues arise when no solution is available at the end of the sampling period. One option to overcome this difficulty is feasible algorithms that guarantee the feasibility of every iterate during the solution process and allow the algorithm to return a suboptimal, but feasible solution.

A successful feasible algorithm is FSQP \cite{Tits1992}, a line-search Sequential Quadratic Programming (SQP) algorithm that solves three subproblems in every iteration to keep the iterates feasible. First, a standard Quadratic Problem (QP) and a modified QP are solved to obtain a feasible descent direction that is close to the standard QP direction, then a second-order correction is performed. An arc search along the combination of both directions yields a feasible direction with a sufficient decrease in the objective function. Notably, the algorithm cannot treat nonlinear equality constraints directly and has to rely on $\ell_1$ relaxed constraints instead \cite{Tits1997}. Additionally, if no feasible initial guess is provided, the algorithm starts with a feasibility phase that attempts to find a feasible point.\\
A recent development is the SEQUOIA solver \cite{Nita2023} that reformulates the general nonlinear optimization problem as a bilevel optimization problem, based on a residual optimization problem. This algorithm also first solves a feasibility problem, which determines an upper bound on the optimal solution. Then, an infeasible, non-tight lower bound of the optimal solution is generated. A bisection search is then performed until a given optimality gap is fulfilled.\\
\cite{Wright2004} introduced the algorithmic framework FP-SQP: it includes a trust-region QP step that is projected onto the feasible set. No particular projection strategy is described, however certain properties should hold at the feasible iterate for the algorithm to globally converge. FP-SQP was successfully used for nonlinear Optimal Control Problems (OCP) with linear constraints in \cite{Tenny2004}, in which a stabilizing feasibility procedure is based on the linear-quadratic regulator gain. The proposed projection works well for the considered OCPs but does not generalize to nonlinear constraints.\\
Another application of FP-SQP was proposed in \cite{Kiessling2022}: a Feasible Sequential Linear Programming (FSLP) algorithm solves trust-region Linear Problems (LPs) and achieves feasible iterates with inner feasibility iterations that are based on iterated second-order corrections. An Anderson acceleration scheme improved the feasibility iteration algorithm in \cite{Kiessling2023}. Several practical deficiencies of the FSLP algorithm were determined during extensive testing. In particular, the requirement to start the algorithm from a feasible point caused difficulties for complex OCPs. Moreover, the LPs in the feasibility iterations could be infeasible, which forces the algorithm to decrease the trust-region radius. Overall, the feasibility iterations mainly converge for small trust-region radii, which produces small steps for the FSLP algorithm.
\vspace{-0.55mm}
\subsection{Contributions}

The contributions of this paper are threefold. Firstly, the practical deficiencies of the FSLP algorithm are demonstrated using illustrative examples. Secondly, we introduce a relaxation of the feasible algorithm based on the tolerance-tube method~\cite{ZoppkeDonaldson1995}, called afSLP. The iterates are not forced to be feasible, but should stay in a tolerance-tube around the feasible set. This allows for early termination of the feasibility iterations, which reduces the overall computational load. Additionally, as in \cite{Tits1992,Nita2023}, a feasibility phase allows for initialization at infeasible points. The phase terminates once an iterate inside of the tolerance-tube is found. Thirdly, we extend the tolerance-tube method with two novel features in order to guarantee global convergence.
Finally, the improved performance of the almost feasible algorithm is demonstrated in practical examples. The tolerance-tube algorithm was implemented in C++ within {CasADi} \cite{Andersson2019}.

\subsection{Outline}

This paper is structured as follows. In Section \ref{sec:FSLP}, the original feasible sequential linear programming algorithm is reviewed. Section \ref{sec:practical_limitations} discusses the practical limitations of FSLP. A relaxed feasibility algorithm based on a tolerance-tube is presented in Section \ref{sec:ToleranceTubeApproach} and simulation results are presented in Section \ref{sec:simulation_results}. Section \ref{sec:conclusion} concludes this paper.

\subsection{Notation}

In order to simplify the presentation, we change the problem formulation compared to \cite{Kiessling2023, Kiessling2022}. The general Nonlinear Problem (NLP) is given by:
\begin{equation}
\label{eq:nlp}
\begin{split}
\min_{w \in \mathbb{R}^{n_w}} \: f(w) \quad \mathrm{s.t.} \:\: g(w)=0,\:\: h(w)\leq 0,
\end{split}
\end{equation}
where $w \in \mathbb{R}^{n_w}$, and $f\colon\mathbb{R}^{n_w}\to\mathbb{R}$, $g\colon\mathbb{R}^{n_w}\to\mathbb{R}^{n_g}$ and $h\colon\mathbb{R}^{n_w}\to\mathbb{R}^{n_h}$ are twice continuously differentiable. The gradient of $f$ at $w$ is given by $\nabla f(w)$. The Jacobian matrices of $g$ and $h$ are represented by $J_g(w) \in \mathbb{R}^{n_g \times n_w}$ and $J_h(w) \in \mathbb{R}^{n_h \times n_w}$, respectively.
The iteration indices are denoted by superscripts, e.g., $w^{(1)}\in\mathbb{R}^{n_w}$ whereas vector components are specified by subscripts, e.g., $w_1 \in \mathbb{R}$. For a given iterate $w^{(k)}$, quantities evaluated at that iterate will be denoted with the same superscript, e.g., $f^{(k)}$ or $\nabla f^{(k)}$.

The measure of infeasibility is defined by $v(w) := \Vert g(w) \Vert_{\infty} + \Vert [h(w)]^+ \Vert_{\infty}$, where $[h(w)]^+ := [\max\{h(w)_i, 0\}]_{i=1}^{n_h}$.
The feasible set is denoted by $\mathcal{F} := \{w \in \mathbb{R}^{n_w} ~|~ g(w) = 0, h(w) \leq 0\}$.

\section{Feasible Sequential Linear Programming}
\label{sec:FSLP}

In this section, FSLP (Algorithm \ref{alg:FPSQP}) is briefly introduced. The feasibility iterations are described in Algorithm \ref{alg:feasibility_iterations}. For a detailed discussion, see \cite{Kiessling2022} and for an extension \cite{Kiessling2023}.

\subsection{Outer Algorithm}

At every iteration, \eqref{eq:nlp} is linearized at the point $w^{(k)}\in\mathcal{F}$. The resulting trust-region LP is:
\begin{equation}
\label{eq:FP-QP}
\begin{aligned}
\min_{w \in \mathbb{R}^{n_w}} \quad &m_f^{(k)}(w):=(\nabla f^{(k)})^{\top} (w-w^{(k)})\\
\mathrm{s.t.}\quad &g^{(k)} + J_g^{(k)} (w-w^{(k)}) = 0,\\
\quad &h^{(k)} + J_h^{(k)} (w-w^{(k)}) \leq 0,\\
\quad &||P(w-w^{(k)})||_{\infty} \leq \Delta^{(k)}.
\end{aligned}
\end{equation}
Here, $P\in\mathbb{R}^{n_y\times n_w}$ denotes a projection and (optionally) scaling matrix that selects the variables involved in the trust region. The identity matrix is usually used.
The solution of \eqref{eq:FP-QP}, $\bar{w}^{(k)}$, is projected onto the feasible set using the feasibility iterations described in Algorithm \ref{alg:feasibility_iterations}. The feasible trial iterate is denoted by $\hat{w}^{(k)}$.
If no feasible point can be found, the trust-region radius is decreased in \eqref{eq:FP-QP}.
The termination criterion is $m_f^{(k)}(\bar{w}) =0$: if no decrease in the model of the objective function $m_f$ at a feasible point is possible, an optimal point was found and the algorithm terminates.

\begin{algorithm}[thpb]
\caption{FSLP}
\label{alg:FPSQP}
\KwIn{Initial point $w^{(0)}\in\mathcal{F}$, projection matrix $P$, initial trust-region radius $\Delta^{(0)} \in (0, \tilde{\Delta}]$, feasibility tolerance $\tau$}
$\mathrm{success} \gets \mathrm{false}$\;
\For{$k=0, 1, 2, \ldots$}{
$\bar{w}^{(k)} \gets$ solve \eqref{eq:FP-QP}\;
\If{$\vert m_f(\bar{w}^{(k)},w^{(k)})\vert \leq \sigma_{\mathrm{outer}}$}{
\textbf{break}
}
$(\hat{w}^{(k)}, \mathrm{success}) \gets \mathrm{FeasIterations}(w^{(k)}, \bar{w}^{(k)}, \Delta^{(k)}, \tau)$ \;
\uIf{$\mathrm{success} = \mathrm{true}$}{
Compute $\rho^{(k)}$ according to \eqref{eq:trustRegionRatioObjective} \;
$\Delta^{(k+1)} \gets \Delta\mathrm{Update}(\rho^{(k)}, w^{(k)}, \bar{w}^{(k)}, \Delta^{(k)})$\;
$(w^{(k+1)}, -) \gets \mathrm{Acceptance}(\rho^{(k)}, w^{(k)}, \hat{w}^{(k)})$
}
\Else{
$\Delta^{(k+1)} \leftarrow \alpha_1 ||P(\bar{w}^{(k)}-w^{(k)})||_{\infty}$\;
$w^{(k+1)} \leftarrow w^{(k)}$\;
}
} 
\Return $w^{(k)}$
\end{algorithm}

The standard trust-region update strategy and the acceptance test are given in Algorithm~\ref{alg:Trust_region_update}.
As usual in trust-region algorithms, the quotient of actual and predicted reduction
\begin{align}
\label{eq:trustRegionRatioObjective}
\rho_{\mathrm{II}}^{(k)} = \Delta f^{(k)} / \Delta m_f^{(k)}  
\end{align}
decides upon step acceptance, where
\begin{equation}
\label{eq:objectiveReductionMetrics}
\Delta f^{(k)} = f(w^{(k)}) - f(\hat{w}^{(k)}),\ 
\Delta m_f^{(k)} = -m_f^{(k)}(\bar{w}^{(k)}).
\end{equation}
The objective serves as merit function since all iterates remain feasible.
Note that $\rho$ has the subscript $II$ for consistency with the extension discussed in Section~\ref{sec:ToleranceTubeApproach}.


\begin{algorithm}[thpb]
\caption{Trust-region radius update}
\label{alg:Trust_region_update}
\SetKwInOut{Parameter}{Parameter}
\Parameter{$\tilde{\Delta} \geq 1$, $0 < \alpha_1 < 1$, $1 < \alpha_2 < \infty$, $0 < \eta_1 < \eta_2 < 1$}
\SetKwProg{myproc}{Procedure}{}{}
\myproc{$\Delta\mathrm{Update}(\rho, w, \bar{w}, \Delta)$}{
    \uIf{$\rho < \eta_1$}{
        \Return $\alpha_1 ||P(\bar{w} - w)||_{\infty}$
    }
    \uElseIf{$\rho > \eta_2$ and $||P(\bar{w} - w)||_{\infty} = \Delta$}{
        \Return $\min(\alpha_2 \Delta, \tilde{\Delta})$
    } 
    \Else{
        \Return $\Delta$
    }
} 
\end{algorithm}
\vspace{-6mm}

\begin{algorithm}[thpb]
\caption{Acceptance test of the trial iterate}
\label{alg:Acceptance_test}
\SetKwInOut{Parameter}{Parameter}
\Parameter{$\sigma \in (0, 1/4)$}
\SetKwProg{myproc}{Procedure}{}{}
\myproc{$\mathrm{Acceptance}(\rho, w, \hat{w})$}{
    \uIf{$\rho > \sigma$}{
        \Return $(\hat{w}, \mathrm{true})$
    }
    \Else{
        \Return $(w, \mathrm{false})$
    }
}
\end{algorithm}


\subsection{Feasibility Iterations}

The feasibility iterations solve a sequence of parametric LPs. In order to distinguish between outer and inner iterations, a second superscript is introduced: $w^{(k,l)}$ denotes the iterate at the $k$-th outer iteration and at the $l$-th inner iteration.

Given $\bar{w}^{(k)}$ the solution of \eqref{eq:FP-QP}, the initial iterate $w^{(k,0)}$ for the feasibility iterations is set to $\bar{w}^{(k)}$ and the parametric linear problem $\mathrm{PLP}(w^{(k,l)}; w^{(k)}, \Delta)$ is defined by:
\begin{equation}
\label{eq:PLP}
\begin{aligned}
\min_{w \in \mathbb{R}^{n_w}} \quad & (\nabla f^{(k)})^{\top} (w - w^{(k)})\\
\mathrm{s.t.} \quad &g(w^{(k,l)}) + J_g^{(k)} (w - w^{(k)}) = 0, \\
\quad & h(w^{(k,l)}) + J_h^{(k)} (w-w^{(k)}) \leq 0, \\
\quad &||P(w - w^{(k)})||_{\infty} \leq \Delta^{(k)}.
\end{aligned}
\end{equation}

Problem \eqref{eq:PLP} is almost identical to \eqref{eq:FP-QP}, except for the constraints evaluated at $w^{(k,l)}$.
Consequently, the feasibility iterations are relatively cheap since only constraint evaluations (no derivatives) are required.

The optimal solution of \eqref{eq:PLP} is $w^{(k,l+1)} := w^*_{\mathrm{PLP}}(w^{(k,l)}; w^{(k)}, \Delta)$. The algorithm continues until an iterate is feasible and fulfills the projection ratio condition, for details see \cite{Kiessling2022}. A heuristic checks whether such an iterate is likely to be found. If this is not the case, the feasibility iterations failed, which results in a decrease of the trust-region radius in Algorithm~\ref{alg:FPSQP}. Only a linear convergence rate is expected from the algorithm.

\begin{algorithm}[thpb]
\caption{Feasibility Iterations}
\label{alg:feasibility_iterations}
\SetKwInOut{Parameter}{Parameter}
\Parameter{$n_{\mathrm{watch}}\in\mathbb{N}$, $\kappa_{\mathrm{watch}}<1$, $\sigma_{\mathrm{inner}}\in (0, 10^{-5})$}
\SetKwProg{myproc}{Procedure}{}{}
\myproc{$\mathrm{FeasIterations}(w^{(k)},\bar{w}^{(k)}, \Delta, \tau)$}{
$w^{(k,0)}\leftarrow \bar{w}^{(k)}$, $\mathrm{success} \leftarrow \mathrm{false}$\;
\For{$l=0, 1, 2, \ldots$}{
\uIf{$v(w^{(k,l)}) \leq \tau$ $\mathrm{and}$ $\Vert \bar{w}^{(k)}-w^{(k,l)}\Vert / \Vert \bar{w}^{(k)}-w^{(k)}\Vert < 1/2$}{
$\mathrm{success} \leftarrow \mathrm{true}$, \textbf{break}\;
}
\ElseIf{$w^{(k,l)}$ diverges according to \cite{Kiessling2022}}{
\textbf{break}\;
}
$w^{(k,l+1)} \gets$ solve $\mathrm{PLP}(w^{(k,l)}, \hat{w}, \Delta)$\;
}
\textbf{return} $(w^{(k, l)}, \mathrm{success})$
} 
\end{algorithm}

\section{Practical Limitations of FSLP}
\label{sec:practical_limitations}

Several practical deficiencies of FSLP, including the feasibility iterations and the initialization at a feasible point, were identified during extensive testing.

\subsection{Feasibility Iterations}
The limitations of the feasibility iterations are illustrated with two examples. The first section addresses the issue of infeasible subproblems due to incompatible constraints. The second problem shows the convergence of the feasibility iterations only for small trust-region radii.
The consequence of these drawbacks is that the FSLP algorithm often takes small steps while guaranteeing the feasibility of every iterate, which causes many additional constraint evaluations and the solution of many additional LPs.

\subsubsection{Infeasible Subproblems}
\begin{figure}[htpb]
\centering
\includegraphics[width=0.7\columnwidth]{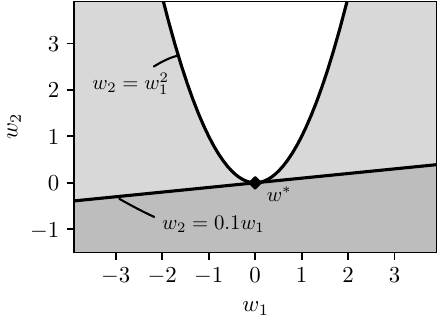}
\caption{Visual representation of Example~\eqref{eq:example_feasibility_iterations}.}
\label{fig:example1}
\end{figure}

Even though the constraints were linearized at a feasible point, the subproblems in subsequent feasibility iterations may be infeasible.
This is illustrated by the test problem:
\begin{align}
\label{eq:example_feasibility_iterations}
\min_{w \in \mathbb{R}^2} w_2 \quad \mathrm{s.t.} \quad w_2 \geq w_1^2, w_2 \geq 0.1 w_1.
\end{align}
Its solution is $w^* = (0, 0)$, as shown in Figure~\ref{fig:example1}. Infeasible areas are grayed out and the solution is shown as a diamond.

\begin{figure}
\centering
\begin{subfigure}{0.5\columnwidth}
\includegraphics[width=\textwidth]{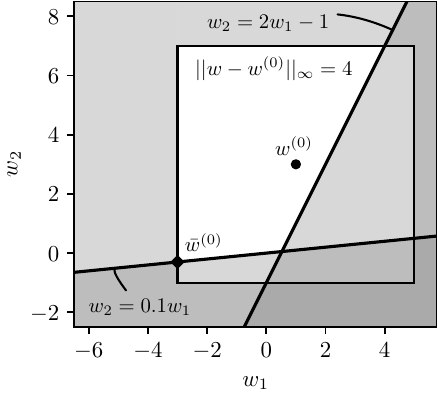}
\caption{First outer iteration.}
\label{fig:example1_feasible_LP}
\end{subfigure}%
\begin{subfigure}{0.5\columnwidth}
\includegraphics[width=\textwidth]{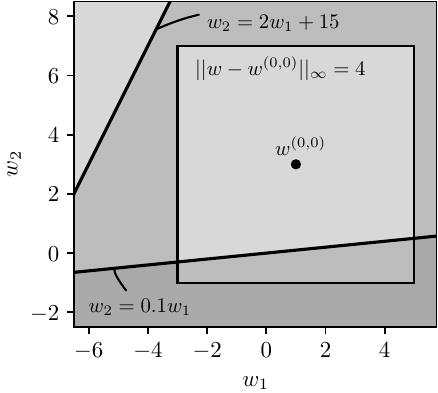}
\caption{First feasibility iteration.}
\label{fig:example1_infeasible_LP}
\end{subfigure}
\caption{Illustration of Feasibility Iterations leading to an infeasible subproblem. Grayed-out areas are infeasible.}
\end{figure}

Since the objective and the second constraint are linear, only the linearization of the first constraint varies over the iterations. Similarly, the trust-region radius
is fixed throughout feasibility iterations.
The FSLP algorithm is started from a point $w^{(0)} = (1, 3)$ with an initial trust-region radius $\Delta^{(0)}=4$. The LP of the first outer iteration is:
\begin{subequations}
\begin{align}
\min_{w \in \mathbb{R}^2} \quad & w_2 \\
\mathrm{s.t.}\quad & w_2 \geq - 1 + 2 w_1, \label{eq:updated_constraint_feasibility_iterations} \\
& w_2 \geq 0.1 w_1, \\
& -3 \leq w_1 \leq 5, -1 \leq  w_2 \leq 7. \label{eq:example_trust_region_constraint} 
\end{align}
\end{subequations}
Its solution is $(-3, -0.3)$, as shown in Figure~\ref{fig:example1_feasible_LP}. Infeasible areas are grayed out and the solution is shown as a diamond.
The linearized constraint of the first feasibility iteration is:
\begin{align}
w_2\geq 15 + 2w_1,
\label{eq:bad_linearized_constraint}
\end{align}
which is a parallel displacement of \eqref{eq:updated_constraint_feasibility_iterations}.
Figure~\ref{fig:example1_infeasible_LP} illustrates that \eqref{eq:bad_linearized_constraint} is shifted outside of the trust region, i.e., the LP becomes infeasible.


We note that this also holds for Anderson Accelerated FSLP since the infeasibility occurs within the first iteration of the feasibility iterations, which is the first iteration for the Anderson acceleration.

\subsubsection{Convergence of Feasibility Iterations for Small Trust-Region Radii}

\begin{figure}[htpb!]
\centering
\includegraphics[width=0.75\columnwidth]{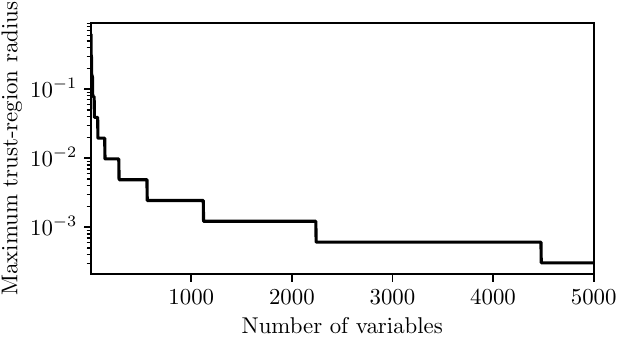}
\caption{Maximum trust-region radii for convergence of feasibility iterations with respect to the number of variables.}
\label{fig:maximum_tr_radii}
\end{figure}

In this section, we analyze the maximum trust-region radius such that the feasibility iterations converge towards a feasible point, while maintaining feasible subproblems and without reaching the maximum number of feasibility iterations.
For $n \geq 2$, we are interested in the following high-dimensional nonlinear problem:
\begin{align*}
\min_{w \in \mathbb{R}^n} \quad -w_1 \quad \mathrm{s.t.} \quad \sum_{i=1}^n w_i^2 - 1 = 0,
\end{align*}
that is finding the maximum $w_1$ on an $n$-dimensional unit sphere.
The feasible initial guess is chosen as $w^{(0)}=(0.5, \sqrt{1-0.5^2}, 0, \ldots, 0)$, i.e., it lies on the two-dimensional unit sphere.
In our implementation, $n \in \{2, 3, 4, \dots, 5000\}$, the feasibility tolerance $\tau$ is set to $10^{-8}$ and the maximum number of feasibility iterations is set to $100$.
The experiment is started with a trust-region radius of $10$ and it is always halved in case the feasibility iterations do not converge.
Figure~\ref{fig:maximum_tr_radii} illustrates the decrease of the maximum trust-region radius that guarantees convergence with the dimension $n$.

This trend was also observed when solving time-optimal control problems (TOCP)
with FSLP. The dimensions of the TOCPs typically range from $100$ to several $1,000$.

\subsection{Feasible Initialization}

In the general case of \eqref{eq:nlp}, it is arduous to initialize the problem at a feasible point. There are two straightforward techniques to obtain a feasible initial guess: i) a feasibility phase that minimizes constraint violation and (hopefully) returns a feasible initial guess~\cite{Tits1997,Nita2023}; or
ii) an $\ell_1$ relaxation of the problem (dropping the equality constraints for simplicity):
\begin{equation*}
\begin{aligned}
\label{eq:l1_relaxed_problem}
\min_{w \in \mathbb{R}^{n_w}} \quad f(w) + \mu e^{\top} s \quad \mathrm{s.t.} \quad g(w) \leq s,\:s \geq 0,
\end{aligned}
\end{equation*}
where $e$ is a vector of ones of appropriate size and $\mu > 0$.
For a given $w$, feasibility is achieved by setting the elastic variables $s$ to sufficiently large values. It is well known that the $\ell_1$ relaxation is exact and that the relaxed problem has the same solution as the original problem \eqref{eq:nlp} provided that the penalty parameter $\mu$ is large enough~\cite{Nocedal2006}.
Moreover, a dynamic penalty parameter update is required for fast convergence. \\
In \cite{Kiessling2023, Kiessling2022}, we decided to relax as few constraints as possible to stay as close as possible to the original problem. In particular, we were interested in solving time-optimal point-to-point motion problems, i.e., initial and terminal conditions are given for the state of the controlled system and the goal is to transfer the system from the initial to the terminal state in minimal time. The considered test problems described obstacle avoidance for one obstacle using an explicit Runge-Kutta 4 integrator for the dynamics. These problems were highly nonlinear and challenging from a mathematical programming perspective, but simple in terms of the system's environment that needed to be controlled. We easily found initial guesses that were feasible with respect to all nonlinear constraints, but not with respect to the initial and terminal conditions. Consequently, the algorithm was started at an almost feasible point. One question remains open: is it worth using the feasibility iterations to keep the problem feasible, even though the iterate is far from the original feasible set? \\
Further testing of FSLP on TOCPs allowed us to identify the following drawbacks with feasible initialization. 
In \cite{Kiessling2022} and \cite{Kiessling2023}, only one obstacle was introduced in the problem. The situation gets more difficult when several obstacles are introduced or if the point-to-point motion is performed in corridors, or in a maze.
If instead of an explicit integrator, we use an implicit integrator, like the implicit Euler method, feasible states are solutions of a system of nonlinear equations, which complicates the search for a feasible initial guess.
The situation becomes even more challenging when a direct collocation approach is chosen, introducing additional collocation points between states. In this case, besides satisfying constraints at states and controls, states are interpolated by a polynomial between two consecutive states, and the dynamics constraints must be met at the collocation point.
This complicates the initialization at a feasible point, leading to the necessity for additional slack variables. In this case, the relaxed problem resembles \eqref{eq:l1_relaxed_problem}, signifying that the relaxed problem is feasible while the original problem is not.

\subsection{Conclusion of Practical Limitations}
From the previous sections, we conclude the following: 
First, the feasibility iterations converge towards a feasible point but only for small trust-region radii, which yields small steps and requires many additional constraint evaluations. In the implementation of FSLP, feasibility is defined with respect to a given tolerance $\tau$ typically set to $10^{-8}$. One straightforward simplification is to relax this tolerance: the iterates are allowed to stay within a tolerance-tube (neighborhood) around the feasible set. This has two advantages: the overall algorithm requires fewer feasibility iterations since feasibility can be achieved only loosely, and the algorithm can take larger steps since it is expected that the feasibility iterations will converge to the required accuracy with a larger trust-region radius.\\
Second, for hard TOCPs, initialization at a feasible point can become difficult and relaxing all the constraints does not seem to be the solution. Additionally, the user experience of FSLP is unfavorable since it is hard to initialize the algorithm. In order to avoid this situation, we propose to initialize the algorithm with a feasibility phase that allows infeasible initial guesses and brings the iterate into a neighborhood around the feasible set. 

\section{An Almost Feasible SLP Algorithm}
\label{sec:ToleranceTubeApproach}

\begin{figure}[htpb!]
\centering
\includegraphics[width=0.8\columnwidth]{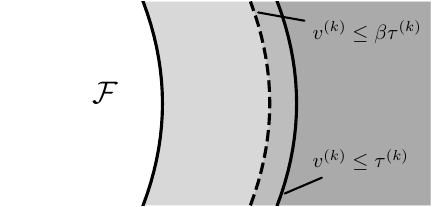}
\caption{Tolerance-tube around the feasible set $\mathcal{F}$. Grayed-out areas are infeasible.}
\label{fig:illustration_tolerance_tube}
\end{figure}

To overcome the issues mentioned in the previous section, we introduce a more general two-phase SLP algorithm, which we call almost feasible Sequential Linear Programming (afSLP), inspired by the tolerance-tube method described in \cite{ZoppkeDonaldson1995}.
A tolerance-tube (Figure~\ref{fig:illustration_tolerance_tube}) is a relaxation of the feasible set $\mathcal{F}$ parameterized by its width $\tau^{(k)} \in (0, 1)$ and a parameter $\beta \in (0, 1)$. $\beta$ is used to guarantee that the tolerance-tube can be reduced while all the iterates stay inside the tolerance-tube.

Phase I improves feasibility until the tolerance-tube is reached. Phase II iterates towards optimality while keeping all iterates inside of the tolerance-tube. When a subproblem is infeasible, a feasibility restoration phase is invoked. Moreover, additional changes were made to the original algorithm that improved the global convergence behavior.

\subsection{Phase I: Feasibility Phase}
Let $w^{(k)}$ be outside of the tolerance-tube, i.e., $v^{(k)} > \beta\tau$.
If the LP~\eqref{eq:FP-QP} is successfully solved, the linearized constraint violation for the solution $\bar{w}^{(k)}$ is zero.
The actual and predicted reductions in constraint violation at $\bar{w}^{(k)}$ are
\begin{equation*}
\Delta v^{(k)} = v^{(k)} - v(\bar{w}^{(k)}),\ \ 
\Delta m_v^{(k)} = v^{(k)}
\end{equation*}
and the ratio of actual and predicted reduction is given by
\begin{equation*}
\rho_{\mathrm{I}}^{(k)}:=\Delta v^{(k)} / \Delta m_v^{(k)}.   
\end{equation*}
The step is accepted or rejected and the trust-region radius is adjusted as described in the FSLP algorithm.
The algorithm stays in Phase I until the iterates reach the tolerance-tube, i.e., $v^{(k)}\leq\beta\tau$, then switches to Phase II.
If however \eqref{eq:FP-QP} is infeasible, the algorithm switches to a feasibility restoration phase described in Section \ref{sec:Restoration_Phase}.

\subsection{Phase II: Optimality Phase}
Let $w^{(k)}$ be inside of the tolerance-tube, i.e., $v^{(k)} \leq \beta \tau$.  If $v(\bar{w}^{(k)}) > \beta \tau$, where $\bar{w}^{(k)}$ is the solution of \eqref{eq:FP-QP}, the feasibility iterations (Algorithm~\ref{alg:feasibility_iterations}) are invoked to bring the trial iterate back in the tolerance-tube. Here, we use an unconstrained sufficient decrease condition (only the objective is considered): the actual and predicted reductions and their quotient are defined as in the FSLP algorithm.

In constrained optimization, the predicted reduction may not be positive, since reducing infeasibility may increase the objective. To address this and to avoid infinitely small steps that converge towards an infeasible point, we use the switching condition traditionally used in filter methods~\cite{Fletcher1999}:
\begin{align}
\label{eq:switchingCondition}
m_f^{(k)}\geq \sigma v^{(k)}
\end{align}
with $\sigma \in (0, 1)$.
If the switching condition is satisfied, the trust-region update strategy (Algorithm \ref{alg:Trust_region_update}) is invoked. Otherwise, sufficient decrease does not seem to be possible, the step is rejected and the trust-region radius is decreased. afSLP terminates if a KKT point is found.
If the subproblem is infeasible, the algorithm switches to the feasibility restoration phase.

\subsection{Feasibility Restoration Phase}
\label{sec:Restoration_Phase}

The feasibility problem is defined as:
\begin{align}
\label{eq:nonsmoothRestoration}
\min_{w \in \mathbb{R}^{n_w}} \quad v_{\mathrm{R}}(w) := \Vert g(w) \Vert_1 + \Vert [h(w)]^+ \Vert_1.
\end{align}
Due to the $\ell_1$ norm, the problem is non-smooth, but by introduction of elastic variables, it can be reformulated as a smooth constrained problem. Let $t^+,t^-\in\mathbb{R}^{n_g}, s\in\mathbb{R}^{n_h}$, then \eqref{eq:nonsmoothRestoration} is equivalent to
\begin{equation}
\label{eq:smoothRestoration}
\begin{aligned}
\min_{w, s, t^+, t^-} \quad & \sum_{i=1}^{n_g} (t^+_i + t^-_i) + \sum_{i=1}^{n_h} s_i \\
\mathrm{s.t.} \quad & g(w) - t^+ + t^- = 0 \\
    & h(w) - s \leq 0 \\
    & s, t^+, t^- \geq 0.
\end{aligned}
\end{equation}
Applying sequential linear programming on \eqref{eq:smoothRestoration}, we obtain the trust-region LP:
\begin{equation}
\label{eq:restoration-LP}
\begin{aligned}
\min_{w, s, t^+, t^-} \quad & \sum_{i=1}^{n_g} (t^+_i + t^-_i) + \sum_{i=1}^{n_h} s_i \\
\mathrm{s.t.} \quad & g^{(k)} + J_g^{(k)} (w - w^{(k)}) - t^+ + t^- = 0,\\
\quad & h^{(k)} + J_h^{(k)} (w - w^{(k)}) - s \leq 0, \\
    & s, t^+, t^- \geq 0, \\
\quad &||P(w-w^{(k)}||_{\infty} \leq \Delta^{(k)}.
\end{aligned}
\end{equation}
Note that elastic variables are excluded from the trust-region constraint so that the LP is feasible.
Let $w_{\mathrm{R}}^{(k)}$ be the $w$ component of the solution of \eqref{eq:restoration-LP}.
We define the model of the $\ell_1$ feasibility measure $v_{\mathrm{R}}$ by 
\begin{align*}
m_{\mathrm{R}}^{(k)}(w) := &\Vert g^{(k)} + J_g^{(k)} (w-w^{(k)}) \Vert_{1}\\
    &+ \Vert [h^{(k)} + J_h^{(k)} (w-w^{(k)})]^+ \Vert_{1},
\end{align*}
the actual and predicted reductions by
\begin{equation*}
\begin{aligned}
\Delta v_{\mathrm{R}}^{(k)} & = v_{\mathrm{R}}^{(k)} - v_{\mathrm{R}}(w_{\mathrm{R}}^{(k)}), \\
\Delta m_{\mathrm{R}}^{(k)} & = v_{\mathrm{R}}^{(k)} - m_{\mathrm{R}}^{(k)}(w_{\mathrm{R}}^{(k)})
\end{aligned}
\end{equation*}
and their ratio by
\begin{equation*}
\rho_{\mathrm{R}}^{(k)} := \Delta v_{\mathrm{R}}^{(k)} / \Delta m_{\mathrm{R}}^{(k)}. 
\end{equation*}
Since \eqref{eq:nonsmoothRestoration} is unconstrained, the predicted reduction is always non-negative and no switching condition is necessary.

The acceptance test (Algorithm \ref{alg:Acceptance_test}) decides upon acceptance or rejection of the trial iterate.
The main difference in our approach compared to \cite{ZoppkeDonaldson1995} is
that $\tau$ is adaptively updated and denoted instead by $\tau^{(k)}$.
If $v^{(k)} < \beta \tau^{(k)}$ and the algorithm switches to feasibility restoration, then we require that $v(\bar{w}^{(k)}) < \beta \tau^{(k)}$. If the trial iterate is accepted, the tolerance-tube width is reduced: $\tau^{(k+1)} \gets \beta \tau^{(k)}$. This is a safeguard that prevents the algorithm from cycling. The parameter $\beta$ guarantees that the tolerance-tube can be reduced while keeping the current iterate inside.

The algorithm terminates when a small step is encountered.
If the current iterate is infeasible, this is an indication that the problem might be infeasible.

\subsection{The Complete Algorithm}

The full afSLP method is described in Algorithm \ref{alg:tolerance_tube}. Since feasibility is not maintained at every iteration, we require for termination that:
\begin{equation*}
\begin{cases}
v(w^*) \leq \varepsilon_{\mathrm{F}} & \text{ (feasibility)} \\
m_f(w^*, \bar{w}^*) \leq \varepsilon_{\mathrm{O}} & \text{ (stationarity)}
\end{cases}
\end{equation*}
for a given feasibility tolerance $\varepsilon_{\mathrm{F}} \in (0, 1)$ and optimality tolerance $\varepsilon_{\mathrm{O}} \in (0, 1)$. If the trust-region radius is decreased below a minimum value, the algorithm terminates with an error message. 

\begin{algorithm}[thpb!]
\caption{afSLP}
\label{alg:tolerance_tube}
\KwIn{Initial point $w^{(0)} \in \mathcal{F}$, projection matrix $P$, initial trust-region radius $\Delta^{(0)} \in (0, \tilde{\Delta}]$, initial tolerance-tube $\tau^{(0)}$, $\beta, \sigma\in (0, 1)$}
\For{$k = 0, 1, 2, \ldots$}{
    $\bar{w}^{(k)} \gets$ solve \eqref{eq:FP-QP}\;
    \uIf{\eqref{eq:FP-QP} is feasible}{
        $\mathrm{update\_radius} \gets \mathrm{false}$ \;
        \uIf(\tcp*[f]{Phase I}){$v^{(k)} > \beta \tau^{(k)}$}{
            $\mathrm{update\_radius} \gets \mathrm{true}$; $\rho^{(k)} \gets \rho^{(k)}_{\mathrm{I}}$
        } 
        \Else(\tcp*[f]{Phase II}){ 
            \If{$v^{(k)} \leq \varepsilon_{\mathrm{F}}$ and $m_f^{(k)}(\bar{w}^{(k)}) \leq \varepsilon_{\mathrm{O}}$}{
                \textbf{break}\;
            }
            \uIf{$v(\bar{w}^{(k)}) \leq \beta \tau^{(k)}$}{
                $\mathrm{success} \gets \mathrm{true}$; $\hat{w}^{(k)} \gets \bar{w}^{(k)}$
            }
            \Else{
                $(\hat{w}^{(k)}, \mathrm{success})\gets\mathrm{FeasIterations}(w^{(k)}, \bar{w}^{(k)}, \Delta^{(k)}, \tau^{(k)})$
            }
            \If{$\mathrm{success}$ and $\Delta m_f^{(k)} \geq \sigma v^{(k)}$}{
                $\mathrm{update\_radius} \gets \mathrm{true}$; $\rho^{(k)} \gets \rho^{(k)}_{\mathrm{II}}$
            }
        }
        \uIf{$\mathrm{update\_radius}$}{
            $\Delta^{(k+1)} \gets \Delta\mathrm{Update}(\rho^{(k)}, w^{(k)}, \bar{w}^{(k)}, \Delta^{(k)})$ \;
            $(w^{(k+1)}, -) \gets \mathrm{Acceptance}(\rho^{(k)}, w^{(k)}, \hat{w}^{(k)})$
        }
        \Else{
            $\Delta^{(k+1)} \leftarrow \alpha_1 ||P(\bar{w}^{(k)} - w^{(k)})||_{\infty}$\;
            $w^{(k+1)} \leftarrow w^{(k)}$\;
        }
    } 
    \Else(\tcp*[f]{Restoration}){
        $\bar{w}^{(k)} \gets$ solve \eqref{eq:restoration-LP}; $\rho^{(k)}\gets\rho^{(k)}_{\mathrm{R}}$\; 
        $\Delta^{(k+1)} \gets \Delta\mathrm{Update}(\rho^{(k)}, w^{(k)}, \bar{w}^{(k)}, \Delta^{(k)})$\;
        $(w^{(k+1)}, \mathrm{accept}) \gets \mathrm{Acceptance}(\rho^{(k)}, w^{(k)}, \bar{w}^{(k)})$\;
        \If{$v^{(k)} \leq \beta \tau^{(k)}$ and $\mathrm{accept}$}{
            $\tau^{(k+1)} \gets \beta\tau^{(k)}$\;
        }
    }
} 
\Return $w^{(k)}$
\end{algorithm}
Global convergence could probably be proved in a similar fashion to \cite{Fletcher1999}. We note that unlike afSLP, the original algorithm does not possess a switching condition~\eqref{eq:switchingCondition}, it does not decrease the tolerance-tube when the feasibility restoration is invoked in Phase II, and it performs a second-order correction, i.e., only one iteration of feasibility iterations, if the trial iterate after solving \eqref{eq:FP-QP} is outside the  tolerance-tube.
We now describe an example where the original tolerance-tube algorithm~\cite{ZoppkeDonaldson1995} would cycle:
\begin{align}
\label{eq:example_cycling_missing_switching_condition}
\min_{w \in \mathbb{R}^2} w_2 \quad \mathrm{s.t.} \quad w_2 \geq w_1^2 + 0.0375, w_1 \geq w_2.
\end{align}
This example is similar to \eqref{eq:example_feasibility_iterations}, except that the quadratic constraint is shifted by $0.0375$ so as to obtain nice iterates in the cycle, and the inequality is flipped in the second constraint while removing the factor 0.1. Without loss of generality, we set the tolerance-tube width to $\tau = 1$ (for lower values of $\tau$, an equivalent problem can be obtained by scaling the constraints -- therefore the constraint violation -- by $\tau$). The initial trust-region radius is set to $\Delta^{(0)} = 1$.

\begin{figure}
\centering
\begin{subfigure}{0.5\columnwidth}
\includegraphics[width=\textwidth]{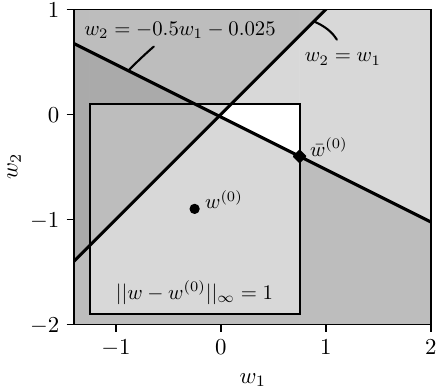}
\caption{First iteration.}
\label{fig:cycling_LP1}
\end{subfigure}%
\begin{subfigure}{0.5\columnwidth}
\includegraphics[width=\textwidth]{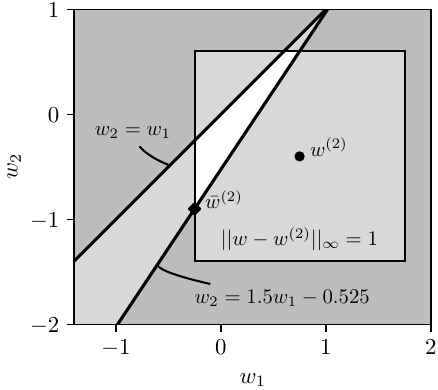}
\caption{Second iteration.}
\label{fig:cycling_LP2}
\end{subfigure}
\caption{Cycling in the original tolerance-tube algorithm.}
\label{fig:tolerance_tube_cycling}
\end{figure}

We start the original tolerance-tube method at $w^{(0)} = (-0.25, -0.9)$. Since $v^{(0)} = 1$, the initial point lies on the tolerance-tube, therefore the algorithm starts in Phase~II. The feasible set of the initial LP is illustrated in Figure~\ref{fig:cycling_LP1} and its solution is $\bar{w}^{(0)} = (0.75, -0.4)$. The trial iterate is accepted and the trust-region radius is increased: $(w^{(1)}, \Delta^{(1)}) = (\bar{w}^{(0)}, 2)$.
The solution of the next LP and the second-order corrected iterate lie outside of the tolerance-tube. Therefore, the trial iterate is rejected and the trust-region radius is decreased: $(w^{(2)}, \Delta^{(2)}) = (w^{(1)}, 1)$. The feasible set of this LP is shown in Figure~\ref{fig:cycling_LP2} and its solution is $\bar{w}^{(2)} = (-0.25, -0.9)$. The trial iterate is accepted and the trust-region radius is again increased: $(w^{(3)}, \Delta^{(3)}) = (\bar{w}^{(2)}, 2)$.
At this point, the original tolerance-tube algorithm starts cycling.

The predicted reduction reveals that $\Delta m_f^{(0)}(\bar{x}^{(0)}) = -0.5$ and $\Delta m_f^{(2)}(\bar{x}^{(0)}) = 0.5$, but $\rho^{(0)} = \rho^{(2)} = 1$: even though the algorithm does not predict a decrease of the objective in the first LP, the trial iterate is accepted since $\rho^{(k)}$ is always $1$. The switching condition~\eqref{eq:switchingCondition} introduced in afSLP prevents the algorithm from taking these steps.


\section{Simulation Results}
\label{sec:simulation_results}

This section presents simulation results that demonstrate the effectiveness of afSLP on two problems: (i) the test problem in \eqref{eq:example_cycling_missing_switching_condition}, and (ii) the time-optimal point-to-point motion (in the Euclidean space) of a parallel SCARA robot.


The OCPs considered in this simulation are implemented using {IMPACT} \cite{Florez2023}, a framework for rapid prototyping of nonlinear model predictive control built on top of {CasADi} which inherits the solvers implemented therein. All simulations were carried out using an Intel core i7-10810U CPU and the LP solver {CPLEX} version 12.8~\cite{cplex2017v12}. The algorithmic parameters are chosen as in \cite{Kiessling2023}.

\subsection{Test problem~\eqref{eq:example_cycling_missing_switching_condition}}
\begin{figure}[htpb]
\centering
\includegraphics[width=\columnwidth]{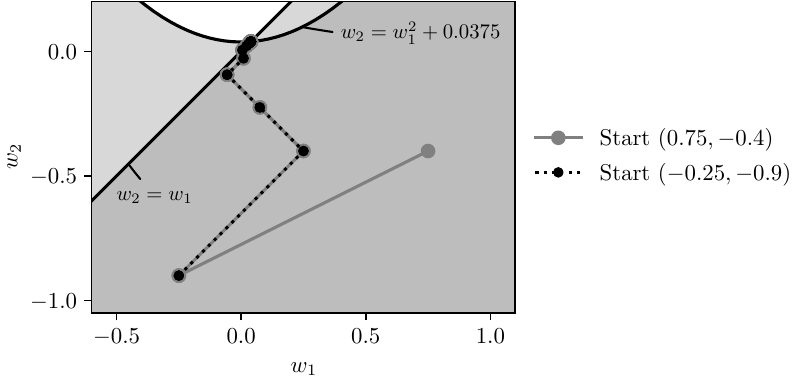}
\caption{Convergence of afSLP algorithm on Problem~\eqref{eq:example_cycling_missing_switching_condition}. Grayed-out areas are infeasible.}
\label{fig:tolerance_tube_cycling_convergence}
\end{figure}

We demonstrate that, unlike the original tolerance-tube algorithm~\cite{ZoppkeDonaldson1995}, afSLP achieves convergence on Problem~\eqref{eq:example_cycling_missing_switching_condition}. 
Since the original tolerance-tube algorithm only checks for $v^{(k)} \leq \tau$ with $\tau = 1$ in Problem~\eqref{eq:example_cycling_missing_switching_condition}, we choose $\tau^{(0)} = 1.2$ and $\beta = 0.9$ to have a similar maximum allowed constraint violation in the initial iteration. The algorithm is started from both cycling points $(-0.25, -0.9)$ and $(0.75, -0.4)$. Figure~\ref{fig:tolerance_tube_cycling_convergence} shows the convergence of afSLP towards the optimal solution $w^* \approx (0.039, 0.039)$ for both initial guesses. As in the analytical example, starting at $(0.75, -0.4)$ brings the next iterate to $(-0.25, -0.9)$, then the iterations for both starting points coincide.

\subsection{Time-optimal motion of a parallel SCARA robot}

To further demonstrate the effectiveness of afSLP, we apply it to the TOCP for the SCARA motion example originally introduced in our previous work~\cite{Kiessling2023}. For the sake of conciseness, only a brief description in the original notation is given below.

The parallel SCARA robot consists of two 2-link arms whose end effectors are attached to a revolute joint $j_5$ (loop closure constraint). It has two actuated joints ($j_1$ and $j_3$) and three unactuated joints ($j_2$, $j_4$ and $j_5$). Let us denote the configuration of the $i$-th joint as $q_i \in \mathbb{S}^1$, the independent coordinates as $q := [q_1,q_3]^\top$, the generalized coordinates as $\bar{q} := \sigma(q) = [q_1,q_2(q),q_3,q_4(q)]^\top$, the inertia matrix as $M$ and the vector of Coriolis and centrifugal effects as ${F}$.
For a state vector $x := [q^{\top},\dot{q}^{\top}]^\top$ and a control (torque) input vector $u := \tau_\mathrm{S} \in \mathbb{R}^2$, the dynamics of the robot~\cite{Cheng_2011} are described by
\begin{equation*}
\dot{x} = f_{\mathrm{ode}}(x,u) := [\dot{q}^\top,\ ({M}^{-1}(\bar{q})(\tau_\mathrm{S}  - {F}(\bar{q},\dot{\bar{q}})))^\top]^\top.
\end{equation*}
It is discretized for a sampling time $h_t$ using a 4th-order Runge-Kutta integrator, which results in a function $f(x_k,u_k,h_{t})$.
The TOCP for the SCARA motion example is formulated as
\begin{subequations}
\label{OCP}
\begin{align}
\min_{\substack{x_0, \ldots, x_N \\ u_0, \ldots, u_{N-1} \\ s_0, s_N, T}} & T + \mu_0^{\top} s_0 + \mu_N^{\top} s_N \\
\mathrm{s.t.} \quad & -s_0 \leq x_0 - \overline{x}_0 \leq s_0, \label{eq:initial_constraint} \\
\quad & x_{k+1} = f(x_k, u_k, h_{t}), & & \hspace{-6.5mm} k \in [0, N-1], \label{eq:ms_constraint} \\
\quad & u_k \in \mathbb{U}, & & \hspace{-6.5mm} k \in [0, N-1], \label{eq:input_constraint} \\
\quad & x_{k} \in \mathbb{X}, & & \hspace{-6.5mm} k \in [0, N], \label{eq:state_constraint} \\
\quad & e(x_k, u_k) \leq 0, & & \hspace{-6.5mm} k \in [0, N-1], \label{eq:stage_constraint} \\
\quad & -s_N \leq x_N - \overline{x}_N \leq s_N \label{eq:terminal_constraint},
\end{align}
\end{subequations}
where $N \in \mathbb{N}$ is the horizon length, $T \in \mathbb{R}_{>0}$ is the time horizon, $h_{t} := T/N$, $x_k \in \mathbb{R}^{n_x}, u_k \in \mathbb{R}^{n_u}$ are the state and control variables, $\bar{x}_0,\bar{x}_N\in\mathbb{R}^{n_x}$ are the initial and terminal states, $\mu_0, \mu_N \in \mathbb{R}^{n_s}_{> 0}$ are penalty parameters for the elastic variables $s_0, s_N \in \mathbb{R}^{n_x}$, and the sets $\mathbb{U}, \mathbb{X}$ are convex polytopes defined by lower and upper bounds on $\tau$, $q$ and $\dot{q}$. Finally, the function $e(x_k, u_k)$ encloses task-related stage constraints: (i) an upper bound $V_{\max}$ on the squared $\ell_2$ norm of the velocity of the end-effector $\dot{p}_\mathrm{ee}$ $(\Vert\dot{p}_\mathrm{ee}\Vert^2 \leq V_{\max}^2)$, and (ii) a collision avoidance constraint based on separating hyperplanes:
\begin{equation}
\mathbf{n}_{\mathrm{a}}^\top p_{\mathrm{ee}} + \mathbf{n}_{\mathrm{b}} + r_{\mathrm{s}} \leq 0,~ \mathbf{n}_{\mathrm{a}}^\top \upsilon + \mathbf{n}_{\mathrm{b}} \geq 0,~ \forall \upsilon \in \mathcal{V}_{\mathrm{obs}},
\end{equation}
where $\mathbf{n} := [\mathbf{n}_{\mathrm{a}}^\top, \mathbf{n}_{\mathrm{b}} ]^\top \in \{\mathbf{v} \in \mathbb{R}^3 : \Vert\mathbf{v}\Vert_{\infty} \leq 1\}$ defines the hyperplane, $r_{\mathrm{s}} \in \mathbb{R}_{\geq 0}$ is a safety margin, and $\mathcal{V}_{\mathrm{obs}}$ represents the vertices of an obstacle.

We compare afSLP against the original FSLP algorithm~\cite{Kiessling2023} and the state-of-the-art nonlinear optimization solver IPOPT~\cite{Waechter2006} for a variation of problem \eqref{OCP} in which the initial and terminal constraints \eqref{eq:initial_constraint} and \eqref{eq:terminal_constraint} are not relaxed.

\subsubsection{Comparison against FSLP}
Table \ref{tab:Mean_eval_g} shows how 
the performance (number of constraint evaluations and average wall time for 100 runs of the problem) of FSLP and afSLP varies for different values of $\tau^{(0)}$. The number of constraint evaluations and the average wall time decrease for increasing values of $\tau^{(0)}$ until $\tau^{(0)} = 10^{-4}$, above which the performance decreases. Overall, afSLP reduces the wall time by 50\%.

\begin{table}[thbp]
\caption{Number of constraint evaluations and average wall time and  on the SCARA test problem.}
\label{tab:Mean_eval_g}
\centering
\begin{tabular}{l c | c c c c c c}
\toprule
 & FSLP & \multicolumn{5}{c}{afSLP} \\
$\tau^{(0)}$ & $10^{-8}$ & $10^{-7}$ & $10^{-6}$ & $10^{-5}$ & $10^{-4}$ & $10^{-3}$\\
\midrule
\# $g, h$ eval. & 723 & 612  & 502 & 557 & 316 & 268 \\  
wall time (s) & 0.576     & 0.447     & 0.361 & 0.423 & 0.283 & 0.287 \\ 
\bottomrule
\end{tabular}
\end{table}

\begin{figure}[htpb!]
\centering
\includegraphics[width=\columnwidth]{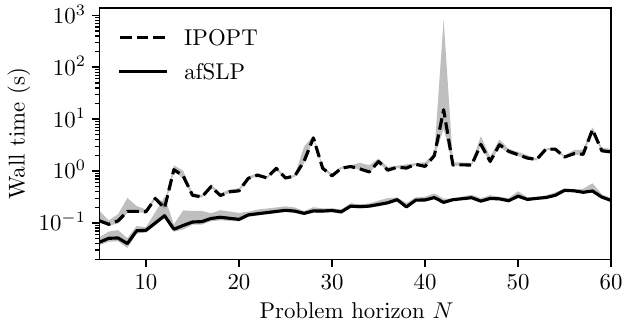}
\caption{Comparison of average wall time of IPOPT and afSLP for different problem horizons $N$.}
\label{fig:cpu_times}
\end{figure}

\subsubsection{Comparison against IPOPT without constraint relaxation}

afSLP is compared against IPOPT on the SCARA test problem for different time horizons $N$.
We choose $\tau^{(0)} = 10^{-3}$, $\beta = 0.9$, $\varepsilon_{\mathrm{O}} = \varepsilon_{\mathrm{F}} = 10^{-7}$. The optimality and feasibility tolerances of IPOPT are set to $10^{-7}$. The maximum iteration number for both solvers was set to $1,000$. For every $N$ the experiment was run $100$ times.
Figure~\ref{fig:cpu_times} shows the average wall times for IPOPT and afSLP. afSLP is significantly faster than IPOPT, which is also due to the fact that expensive evaluations of the Hessian matrix are not required.

\section{Conclusion}
\label{sec:conclusion}

In this paper, we proposed illustrative examples that exhibit the practical deficiencies of FSLP, a feasible sequential linear programming method for solving TOCPs. We then introduced afSLP, an almost feasible SLP method that alleviates the FSLP shortcomings: it allows infeasible initial guesses, iterates first towards the feasible set until a tolerance-tube is reached, then iterates towards optimality. It is inspired by a tolerance-tube method introduced in \cite{ZoppkeDonaldson1995}. 
We demonstrated in a simple example the lack of global convergence of the original tolerance-tube method and proposed two enhancements that overcome these limitations. The performance of afSLP was assessed on a SCARA robot optimization problem and compared against FSLP and the state-of-the-art nonlinear optimization solver IPOPT. \\
Future work will focus on establishing a global convergence proof, on improving the local convergence behavior and on providing suboptimal feasibility guarantees even though afSLP maintains a relaxed feasible set.

\bibliography{bibliography}

\begin{thebibliography}{10}

\bibitem{Andersson2019}
J.~A.~E. Andersson, J.~Gillis, G.~Horn, J.~B. Rawlings, and M.~Diehl.
\newblock {CasADi} -- a software framework for nonlinear optimization and
  optimal control.
\newblock {\em Math. Program. Comput.}, 11(1):1--36, 2019.

\bibitem{Cheng_2011}
L.~Cheng, Y.~Lin, Z.-G. Hou, M.~Tan, J.~Huang, and W.~J. Zhang.
\newblock Adaptive tracking control of hybrid machines: A closed-chain five-bar
  mechanism case.
\newblock {\em {IEEE}/{ASME} Trans. Mechatron.}, 16(6):1155--1163, dec 2011.

\bibitem{Fletcher1999}
R.~Fletcher, S.~Leyffer, and P.~Toint.
\newblock On the global convergence of an {SLP}-filter algorithm.
\newblock {\em Numerical Analysis Report NA/183, University of Dundee, UK}, 12
  1999.

\bibitem{Florez2023}
A.~Florez, A.~Astudillo, W.~Decré, J.~Swevers, and J.~Gillis.
\newblock {IMPACT}: A toolchain for nonlinear model predictive control
  specification, prototyping, and deployment.
\newblock In {\em Proceedings of the IFAC World Congress}, Yokohama, Japan,
  2023.

\bibitem{cplex2017v12}
Cplex~IBM ILOG.
\newblock V12.8: User’s manual for {CPLEX}, 2017.

\bibitem{Kiessling2023}
D.~Kiessling, P.~Pas, A.~Astudillo, P.~Patrinos, and J.~Swevers.
\newblock Anderson accelerated feasible sequential linear programming.
\newblock In {\em Proceedings of the IFAC World Congress}, Yokohama, Japan,
  2023.

\bibitem{Kiessling2022}
D.~Kiessling, A.~Zanelli, A.~Nurkanovi\'c, J.~Gillis, M.~Diehl, M.~Zeilinger,
  G.~Pipeleers, and J.~Swevers.
\newblock A feasible sequential linear programming algorithm with application
  to time-optimal path planning problems.
\newblock In {\em Proceedings of 61st IEEE Conference on Decision and Control},
  Cancun, Mexico, December 2022.

\bibitem{Nita2023}
L.~Nita and E.~Carrigan.
\newblock {SEQUOIA}: A sequential algorithm providing feasibility guarantees
  for constrained optimization.
\newblock In {\em Proceedings of the IFAC World Congress}, Yokohama, Japan,
  2023.

\bibitem{Nocedal2006}
J.~Nocedal and S.~J. Wright.
\newblock {\em Numerical Optimization}.
\newblock Springer Series in Operations Research and Financial Engineering.
  Springer, 2 edition, 2006.

\bibitem{Tenny2004}
M.~J. Tenny, S.~J. Wright, and J.B. Rawlings.
\newblock {N}onlinear model predictive control via feasibility-perturbed
  sequential quadratic programming.
\newblock {\em Comput. Optim. Appl.}, 28:87--121, 2004.

\bibitem{Tits1997}
A.~Tits and C.~Lawrence.
\newblock Nonlinear equality constraints in feasible sequential quadratic
  programming.
\newblock {\em Optimization Methods and Software}, 6, 02 1997.

\bibitem{Waechter2006}
A.~W\"achter and L.~T. Biegler.
\newblock On the implementation of an interior-point filter line-search
  algorithm for large-scale nonlinear programming.
\newblock {\em Mathematical Programming}, 106(1):25--57, 2006.

\bibitem{Wright2004}
S.~J. Wright and M.~J. Tenny.
\newblock {A} feasible trust-region sequential quadratic programming algorithm.
\newblock {\em SIAM Journal on Optimization}, 14:1074--1105, 1 2004.

\bibitem{Tits1992}
J.~L. Zhou and A.~Tits.
\newblock User's guide for {FSQP} version 3.0c: A {FORTRAN} code for solving
  constrained nonlinear (minimax) optimization problems, generating iterates
  satisfying all inequality and linear constraints.
\newblock 01 1992.

\bibitem{ZoppkeDonaldson1995}
C.~Zoppke-Donaldson.
\newblock {\em A tolerance-tube approach to sequential quadratic programming
  with applications}.
\newblock PhD thesis, University of Dundee, 1995.

\end{thebibliography}
\bibliographystyle{plain} 
\addtolength{\textheight}{-12cm}   







\end{document}